\documentstyle[12pt]{article}
\newcommand{\const}{\mathop{\rm const}\limits}

\newcommand{\supp}{\mathop{\rm supp}\limits}

\newcommand{\diam}{\mathop{\rm diam}\limits}

\newcommand{\mes}{\mathop{\rm mes}\limits}

\newcommand{\card}{\mathop{\rm card}\limits}

\newcommand{\cov}{\mathop{\rm cov}\limits}

\begin{document}

\begin{center}

{\bf  CENTRAL LIMIT THEOREM IN H\"OLDER SPACES } \\

\vspace{3mm}

{\bf IN THE TERMS OF  MAJORIZING MEASURES.} \par

\vspace{4mm}

 $ {\bf E.Ostrovsky^a, \ \ L.Sirota^b } $ \\

\vspace{4mm}

$ ^a $ Corresponding Author. Department of Mathematics and computer science, Bar-Ilan University, 84105, Ramat Gan, Israel.\\
 E-mail:  \  eugostrovsky@list.ru\\

\vspace{4mm}

$ ^b $  Department of Mathematics and computer science. Bar-Ilan University,
84105, Ramat Gan, Israel.\\

E-mail: \ sirota3@bezeqint.net\\

\vspace{7mm}

                    {\sc Abstract.}\\

 \end{center}

 \vspace{3mm}

 We obtain some sufficient conditions  for the Central Limit Theorem for the random processes (fields) with values
in the separable part of  H\"older  space in the modern terms of  majorizing (minorizing) measures,
belonging to X.Fernique and M.Talagrand.\par

 We introduce a new class of Banach spaces-rectangle H\"older spaces and investigate CLT
 in this spaces via the fractional order Sobolev-Grand Lebesgue norms. \par

 Our further considerations based on the improvement of the L.Arnold and P.Imkeller generalization of the classical
Garsia-Rodemich-Rumsey inequality, which allow us to reduce degree of the distance in the important particular cases. \par

  \vspace{4mm}

\par

\vspace{4mm}

 \par

\vspace{4mm}

\section{Notations. H\"older spaces. Statement of problem. History.}

\vspace{3mm}

 Let $ (X = \{x \},d) $ be compact metric space relative some distance (or semi - distance) $  d = d(x_1, x_2). $
 The  H\"older (Lipshitz) space  $ H^o(d)  $  consists by definition on all the numerical (real or complex)
 continuous relative the distance $ d = d(t,s) $ functions $ f: T \to R $  satisfying the condition

$$
\lim_{\delta \to 0+} \frac{\omega(f,d, \delta)}{\delta} = 0.  \eqno(1.1)
$$
 Here $ \omega(f, \delta)  $ is uniform module of continuity of the (continuous) function $  f:$

$$
\omega(f,d, \delta) = \omega(f,\delta) = \sup_{t,s: d(t,s) \le \delta} |f(t) - f(s)|. \eqno(1.2)
$$

The norm of the space $ H^o(\omega) $  is defined as follows:

$$
||f||H^o(d) = \sup_{t \in T} |f(t)| + \sup_{d(x_1, x_2) > 0} \left\{ \frac{|f(x_1) - f(x_2)|}{d(x_1, x_2)}  \right\}. \eqno(1.3)
$$

 The detail investigation of these spaces with applications in the theory of non - linear singular integral equations
  is undergoing in the first chapter of a monograph  of  Gusejnov A.I., Muchtarov Ch.Sh.
   \cite{Gusejnov1}. We itemize some used facts about these spaces.\par

 This modification of the classical H\"older (Lipshitz) space is Banach space, i.e. is  linear, normed, complete and separable. \par

 Note but the space $ H^o(d) $  may be trivial, i.e. may consists only constant functions. Let for instance, $  X  $ be convex
connected closed bounded domain in the space $  R^m, \ m = 1,2, \ldots $ and let $  d(x_1, x_2) $ be usual Euclidean distance.
Then the space $ H^o(d)  $ is trivial: $  \dim H^o(d) = 1. $\par
  The space $  H^o(d^{\beta}), \ \beta = \const \in (0,1) $  in this example in contradiction is  not trivial. \par

 \vspace{3mm}

Further, if an another distance $ r = r(x_1,x_2)  $ on the source set  $ X $ is such that

$$
\forall x_1 \in X \ \Rightarrow  \lim_{d(x,x_1) \to 0} \frac{d(x,x_1)}{r(x,x_1)} = 0, \eqno(1.4)
$$
then the space $ H^o(d) $ is continuously  embedded in the space $ H^o(r). $ \par
 We will write the equality (1.4) as follows: $  d << r. $\par

 For instance, the distance $  r(x_1,x_2) $ may has a form

$$
r(x_1, x_2) = d^{\beta}(x_1, x_2),  \ \beta = \const \in (0,1).
$$

\vspace{3mm}

 Let $  \xi = \xi(x), \ x \in X $ be {\it in the sequel, during whole article } be
 separable numerical centered (mean zero) random process (r.pr) or equally random field (r.f.)
with finite (bounded) covariation function

$$
R(x_1, x_2) = R_{\xi}(x_1, x_2) = \cov (\xi(x_1), \xi(x_2)) = {\bf E} \xi(x_1) \cdot \xi(x_2). \eqno(1.5)
$$
 Let also $  \xi_i(x), \ i=1,2,3,\ldots $ be independent copies of r.f. $  \xi(x), $ defined may be on some
sufficiently rich probability space,

$$
S_n(x) = n^{-1/2} \sum_{i=1}^n \xi_i(x).
$$

\vspace{3mm}

 Evidently, the finite - dimensional distributions of the sequence of the r.f. $ S_n(x) $ converge as $  n \to \infty $ to the
finite - dimensional distribution of the Gaussian separable mean zero r.f. $ S_{\infty}(x) $  with at the same covariation function
$ R_{\xi}(x_1, x_2). $\\

\vspace{3mm}

{\bf Definition 1.1.}  The r.f. $ \xi(x) $ or equally the sequence of normed r.f.  $ \{  \xi_i(x) \} $ satisfies by definition
the Central Limit Theorem (CLT) in the space $  H^o(d)  $ (or analogously in arbitrary another separable Banach space) iff

\vspace{3mm}

{\bf 1.}   $ {\bf P} \left(\xi(\cdot) \in H^o(d) \right) = 1; $ \\

\vspace{3mm}

{\bf 2.}  The limiting Gaussian r.f. $  S_{\infty}(x) $ belongs also to this apace $ H^o(d) $ a.e. \\

\vspace{3mm}

{\bf 3.}  The sequence of distributions of the r.f. $ S_n(\cdot) $ in the space $  H^o(d) $ converges weakly as $  n \to \infty $
to the distribution of the r.f. $ S_{\infty}(x). $ \\

 \vspace{3mm}

 The last statement denotes that for arbitrary continuous bounded functional $ F: H^o(d) \to R $

$$
\lim_{n \to \infty} {\bf E} F(S_n(\cdot)) = {\bf E} F(S_{\infty}).
$$

 In particular,

$$
\lim_{n \to \infty} {\bf P} ( ||S_n(\cdot)||H^o(d) > u) = {\bf  P} (||S_{\infty}||H^o(d) > u), \ u > 0.
$$

\vspace{4mm}

{\bf  Our aim in this article is obtaining some sufficient condition for CLT in H\"older space in the too modern terms of
 majorizing (minorizing) measures. } \par

\vspace{3mm}

 There are many works  containing the CLT in Banach  spaces, see e.g.  monographs \cite{Dudley1}, \cite{Ledoux1},
\cite{Ostrovsky1}.  The recent version for CLT in H\"older spaces, for example for the Banach space valued random processes,
formulated in fact in the entropy terms see in  \cite{Ostrovsky1}, chapter 4, section 4.13. (1999);
 \cite{Ratchkauskas1}-\cite{Ratchkauskas5} (2004-2006); \cite{Klicnarov'a1}, (2007). \par

 In the article \cite{Ratchkauskas2} is obtained the necessary and sufficient condition in entropy terms
  for the H\"olderian functional central limit theorem. \par

 A very important applications of this CLT  in the   epidemic change statistics is described in
\cite{Ratchkauskas3}, \cite{Ratchkauskas4}. Another possible applications for the functional CLT appears in the parametric Monte-Carlo method,
\cite{Frolov1}, \cite{Grigorjeva1}, \cite{Ostrovsky106}.  \par

 In the article of  B.Heinkel \cite{Heinkel1}  is obtained sufficient condition for  CLT in the space of continuous functions
$  C(T,d) $  in the more modern and more strong terms of "majorizing measures" or equally "generic chaning";  see \cite{Dudley1},
\cite{Fernique1}, \cite{Talagrand1}- \cite{Talagrand4}. \par

 Notice that  the CLT in the space $  C(T,d) $ follows the CLT in {\it some } H\"older space $ H^o(r), \ d << r. $
\cite{Ostrovsky1}, chapter 4, section 4.13. \par

 It is interest by our opinion to obtain the conditions for CLT also in the H\"older spaces  in these terms. \par

\vspace{4mm}

\section{ Majorizing and minorizing measures.}

\vspace{4mm}

 We recall here for reader convenience some used further facts about the theory of majorizing and minorizing measures.
This classical definition with theory explanation and applications basically in the investigation of local structure
of random processes and fields  belongs to X.Fernique \cite{Fernique1},  \cite{Fernique2},  \cite{Fernique3} and M.Talagrand
\cite{Talagrand1}, \cite{Talagrand2}, \cite{Talagrand3}, \cite{Talagrand4}, \cite{Talagrand5}.
 See also \cite{Bednorz1}, \cite{Bednorz2}, \cite{Bednorz3}, \cite{Dudley1}, \cite{Ledoux1},
 \cite{Ostrovsky100}, \cite{Ostrovsky101}, \cite{Ostrovsky102}, \cite{Ostrovsky107}. \par

 \vspace{3mm}

 Let $ (X,d),  (Y,\rho) $ be separable metric spaces, $ m  $ be arbitrary distribution, i.e.
 Radon probabilistic measure on the  set $ X, $
$ f: X \to Y $  be (measurable) function. Let also $ \Phi(z), \ z \ge 0 $ be continuous Young-Orlicz function,  i.e.
strictly increasing function such that

$$
\Phi(z) = 0 \ \Leftrightarrow z = 0; \ \lim_{z \to \infty} \Phi(z) = \infty.
$$
 We denote as usually

$$
\Phi^{-1}(w) = \sup \{z, z \ge 0, \  \Phi(z) \le w \}, \ w \ge 0
$$
the inverse function to the function $  \Phi; $
$$
B(d,r,x) =  B(r,x) = \{ x_1: \ x_1 \in X, \ d(x_1,x) \le r  \}, \ x \in X, \ 0 \le r \le \diam(X)
$$
be the closed ball of radii $ r $ with center at the point $ x. $\par
 Let us introduce the Orlicz space $  L(\Phi) = L(\Phi; m \times m, \ X \otimes X)  $ on the set $ X \otimes X $ equipped with
the Young-Orlicz function $ \Phi. $ \par

 \vspace{4mm}

{\it We assume henceforth that for all the values } $ x_1, x_2 \in X, \ x_1 \ne x_2 $ {\it (the case } $ x_1 = x_2 $  {\it is trivial)}
{\it the value  } $ \rho(f(x_1), f(x_2)) $ {\it belongs to the space } $ L(\Phi). $ \par

 As a rule,

$$
\rho(f(x_1), f(x_2))  = |f(x_1) - f(x_2)|.
$$

\vspace{4mm}

 Note that for the existence of such a function $ \Phi(\cdot) $ is necessary and sufficient only the integrability
of the distance   $ \rho(f(x_1), f(x_2)) $ over the product measure $ m \times m:$

$$
\int_X \int_X \rho(f(x_1), f(x_2)) \ m(dx_1) \  m(dx_2) < \infty,
$$
see \cite{Krasnoselsky1}, chapter 2, section 8. \par

  Under this assumption the distance $ d = d(x_1, x_2) $ may be constructively defined by the formula:

$$
d_{\Phi}(x_1,x_2) := || \rho(f(x_1), f(x_2))||L(\Phi), \eqno(2.1)
$$
where $  || \cdot ||L(\Phi) $ denotes the  Orlicz's norm. \par

  Since the function $ \Phi = \Phi(z) $ is presumed to be continuous and strictly  increasing,
it follows from the relation (1.1) that $ V(d_{\Phi}) \le 1, $ where by definition

$$
V(d):= \int_X \int_X \Phi \left[ \frac{\rho(f(x_1), f(x_2))}{d(x_1,x_2)} \right] \ m(dx_1) \  m(dx_2). \eqno(2.2)
$$

 Let us define also the following important distance function: $ w(x_1, x_2) = $
 $$
  w(x_1, x_2; V) = w(x_1, x_2; V, m ) = w(x_1, x_2; V, m,\Phi) = w(x_1, x_2; V, m,\Phi,d) \stackrel{def}{=}
 $$

$$
 6 \int_0^{d(x_1, x_2)} \left\{ \Phi^{-1} \left[ \frac{4V}{m^2(B(r,x_1))} \right] +
\Phi^{-1}  \left[ \frac{4V}{m^2(B(r,x_2))} \right] \right\} \ dr, \eqno(2.3)
$$
 where $ m(\cdot) $ is probabilistic Borelian measure on the set $ X.$ \par
  The triangle inequality and other properties of the distance function $ w = w(x_1, x_2) $ are proved in
\cite{Kwapien1}.\par

\vspace{3mm}

{\bf Definition 2.1. } (See  \cite{Kwapien1}). The measure $ m  $ is said to be
{\it minorizing measure } relative the distance $ d = d(x_1,x_2), $ if for each values $ x_1, x_2 \in X
\ V(d) < \infty $ and moreover $ \ w(x_1,x_2; V(d)) < \infty. $\par

\vspace{3mm}

We will denote the set of all minorizing measures on the metric set $ (X,d) $  by  $ \cal{M} = \cal{M}(X).$ \par

 Evidently, if the function $ w(x_1, x_2) $ is bounded, then the minorizing measure $  m  $ is majorizing.  Inverse
proposition is not true, see  \cite{Kwapien1}, \cite{Arnold1}. \par

\vspace{3mm}

{\bf Remark 2.1.} If the measure $  m  $ is minorizing, then

$$
w(x_n, x ; V(d)) \to 0 \ \Leftrightarrow d(x_n, x) \to 0, \ n \to \infty.
$$
 Therefore, the continuity of a function relative the distance $  d  $ is equivalent to
the continuity  of this function  relative the distance $  w.  $ \par

\vspace{3mm}

{\bf Remark 2.2.}  If

$$
\sup_{x_1, x_2 \in X} w(x_1, x_2; V(d)) < \infty,
$$
then the measure $ m $ is called {\it majorizing measure.} \par

 Some considerations about the choice of the majorizing (minorizing) measures see in the article
 \cite{Ostrovsky602}; see also reference therein.\par

\vspace{4mm}

 The following important inequality belongs to  L.Arnold and P.Imkeller \cite{Arnold1}, \cite{Imkeller1};
see also \cite{Kassman1}, \cite{Barlow1}. \par

\vspace{4mm}

{\bf  Theorem of  L.Arnold and P.Imkeller.} {\it Let the measure $  m  $ be minorizing. Then there exists a
modification of the function $ f $ on the set of zero measure, which we denote also by $ f, $ for which}

$$
\rho(f(x_1), f(x_2)) \le  w(x_1, x_2; V, m,\Phi,d). \eqno(2.4)
$$
{\it  As a consequence: this function $  f  $ is $  d - $ continuous and moreover $ w - $ Lipshitz continuous
with unit constant. }\par

\vspace{4mm}

 The inequality (2.4) of L.Arnold and P.Imkeller is significant generalization of celebrated
Garsia-Rodemich-Rumsey inequality, see \cite{Garsia1}, with at the same applications as mentioned before
\cite{Hu1}, \cite{Ostrovsky100}, \cite{Ostrovsky101}, \cite{Ostrovsky102}, \cite{Ral'chenko1}. \par

\vspace{3mm}

{\bf Remark 2.3.} The inequality of  L.Arnold and P.Imkeller (2.4) is closely related with the theory of
fractional order Sobolev's - rearrangement invariant spaces, see  \cite{Barlow1}, \cite{Garsia1}, \cite{Hu1}, \cite{Kassman1},
\cite{Nezzaa1}, \cite{Ostrovsky101}, \cite{Ral'chenko1}, \cite{Runst1}.\par

\vspace{3mm}

{\bf Remark 2.4.} In the previous articles \cite{Kwapien1}, \cite{Bednorz4} was imposed on the function $ \Phi(\cdot) $
the following $ \Delta^2 $ condition:

$$
\Phi(x) \Phi(y) \le \Phi(K(x+y)), \ \exists K = \const \in (1,\infty), \ x,y \ge 0
$$
or equally

$$
\sup_{x,y > 0} \left[ \frac{\Phi^{-1}(xy)}{\Phi^{-1}(x) + \Phi^{-1}(y)}\right] < \infty. \eqno(2.5)
$$
{\it  We do not suppose this condition. For instance, we can consider the function of a view $ \Phi(z) = |z|^p,  $
which does not satisfy the condition (2.5)}. \par

\vspace{4mm}

\section{ H\"older's  CLT over Lebesgue-Riesz spaces. }

\vspace{4mm}

  Let $ \xi = \xi(x), \ x \in X $ be again separable centered continuous {\it in probability }
  random field (r.f),  not necessary to be Gaussian.   The correspondent
   probability and expectation we will denote  by $ {\bf P}, \ {\bf E,}  $  and  the probabilistic
 Lebesgue-Riesz $ L_p $   norm of a random variable (r.v) $  \eta $  we will denote as follows:

 $$
 |\eta|_p \stackrel{def}{=} \left[ {\bf E} |\eta|^p \right]^{1/p}.
 $$

\vspace{3mm}

 Let the r.f. $ \xi(\cdot) $ be such that

 $$
\exists p = \const \ge 2 \ \Rightarrow    \sup_{x \in X} |\xi(x)|_p < \infty.
 $$
Then we can define a so-called natural, or Pisier's distance \cite{Pisier1} $ d_p = d_p(x_1, x_2) $ as follows

$$
d_p(x_1,x_2) \stackrel{def}{=} |\xi(x_1) - \xi(x_2)|_p, \eqno(3.0)
$$
which is evidently bounded.\par

\vspace{3mm}

{\bf Theorem 3.1.}   Suppose  the measure $ m $ and distance $  d_p  $ are such that

$$
 m^2(B(d_p,r,x)) \ge r^{\theta}/C(\theta), \ r \in [0, \diam(X, d_p)], \ \exists \ \theta = \const > 0, \ C(\theta) \in (0,\infty). \eqno(3.1)
$$

 Let also $ p = \const  > \theta, $ so that $  p > \max(\theta, 2). $ \par

 Our statement: for arbitrary (semi -) distance $ \rho = \rho(x_1, x_2) $ such that $ d_p << \rho  $ the r.f.
$ \xi(x) $ satisfies the CLT in H\"older space $  H^o(\rho). $ \par

\vspace{3mm}

{\bf Proof.} We will use the following proposition from the article \cite{Ostrovsky107} (Proposition 2.1.):
we get using the inference also  theorem 2.1 therein that for the r.f. $ \xi = \xi(x) $ the following inequality holds:
$ m \in  \cal{M} $  and

$$
|\xi(x_1) - \xi(x_2)|  \le 12 \ Z^{1/p} \ 4^{1/p} \ C^{1/p}(\theta) \  \frac{d_p^{1-\theta/p}(x_1, x_2)}{1-\theta/p}, \eqno(3.2)
$$
where the r.v. $  Z  $ has unit expectation:  $ {\bf E} Z = 1. $\par
 We intent to apply the inequality (3.2) for the random fields $ S_n(\cdot)  $ instead $ \xi(x). $ Note first of all that
the classical Rosenthal's inequality \cite{Rosenthal1} asserts in particular that if $  \{ \zeta_i \}, \ i = 1,2,\ldots $ are
the sequence of i., i.d. {\it centered} r.v. with finite $  p^{th} $ moment, then

$$
\sup_n \left| n^{-1/2} \sum_{i=1}^n \zeta_i  \right|_p  \le \frac{C_R \ p}{ e \cdot \ln p} \ |\zeta_1|_p, \ p \ge 2. \eqno(3.3)
$$
 About the exact value of the constant $  C_R $ see the article \cite{Ostrovsky601}.  Note that for the symmetrical distributed r.v.
 $  C_R \le 1.53573. $\par

 We have using Rosenthal's inequality since $ p \ge 2 $

$$
| S_n(x_1) - S_n(x_2)|_p \le d_p(x_1, x_2) \cdot \frac{C_R \ p}{e \cdot \ln p} \le C_1(p) \cdot d_p(x_1, x_2),
$$
and we conclude by means of estimate (3.2)

$$
|S_n(x_1) - S_n(x_2)|  \le C_2(\theta,p) \cdot Z_n \cdot  \frac{d_p^{1-\theta/p}(x_1, x_2)}{1-\theta/p} =
 C_3(\theta,p) \cdot Z_n^{1/p} \cdot  d_p^{1-\theta/p}(x_1, x_2), \eqno(3.4)
$$
where $  Z_n $ is the sequence of non - negative  r.v. with unit expectation $ {\bf E} Z_n = 1. $\par
 Let $  \nu = \nu(x_1, x_2) $ be arbitrary  {\it  intermediate  } distance on the set $  X $ between $ r(\cdot, \cdot) $ and
 $ d^{1 - \theta/p}_p(\cdot, \cdot): $

$$
  d^{1 - \theta/p}_p(\cdot, \cdot) <<   \nu(\cdot, \cdot)  <<   r(\cdot, \cdot).
$$
 We deduce from (3.4)

$$
\frac{|S_n(x_1) - S_n(x_2)|}{\nu(x_1,x_2)} \le  C_3(\theta,p) \cdot Z_n^{1/p} \cdot
 \frac{d_p^{1-\theta/p}(x_1, x_2)}{ \nu(x_1,x_2)  }, \eqno(3.5)
$$
and we conclude taking into account the structure of compact embedded  H\"older subspaces into ones that
the sequence of r.f. $ S_n(\cdot) $ satisfies of the famous Prokhorov's criterion \cite{Prokhorov1}
for the weak compactness  of its distributions in the H\"older space $ H^o(r). $\par

 This completes the proof of theorem 3.1. \par

\vspace{3mm}

\section{ Main result: Grand Lebesgue spaces approach. }

\vspace{3mm}

 We recall first of all briefly  the definition and some simple properties of the so-called Grand Lebesgue spaces;   more detail
investigation of these spaces see in \cite{Fiorenza3}, \cite{Iwaniec2}, \cite{Kozachenko1}, \cite{Liflyand1}, \cite{Ostrovsky1},
\cite{Ostrovsky2}; see also reference therein.\par

  Recently  appear the so-called Grand Lebesgue Spaces $ GLS = G(\psi) =G\psi =
 G(\psi; A,B), \ A,B = \const, A \ge 1, A < B \le \infty, $ spaces consisting
 on all  the random variables (measurable functions) $ f: \Omega \to R $ with finite norms

$$
   ||f||G(\psi) \stackrel{def}{=} \sup_{p \in (A,B)} \left[ |f|_p /\psi(p) \right]. \eqno(4.1)
$$

  Here $ \psi(\cdot) $ is some continuous positive on the {\it open} interval
$ (A,B) $ function such that

$$
     \inf_{p \in (A,B)} \psi(p) > 0, \ \psi(p) = \infty, \ p \notin (A,B).
$$
 We will denote
$$
 \supp (\psi) \stackrel{def}{=} (A,B) = \{p: \psi(p) < \infty, \}
$$

The set of all $ \psi $  functions with support $ \supp (\psi)= (A,B) $ will be
denoted by $ \Psi(A,B). $ \par
  This spaces are rearrangement invariant, see \cite{Bennet1}, and
    are used, for example, in the theory of probability  \cite{Kozachenko1},
  \cite{Ostrovsky1}, \cite{Ostrovsky2}; theory of Partial Differential Equations \cite{Fiorenza3},
  \cite{Iwaniec2};  functional analysis \cite{Fiorenza3}, \cite{Iwaniec2},  \cite{Liflyand1},
  \cite{Ostrovsky2}; theory of Fourier series, theory of martingales, mathematical statistics,
  theory of approximation etc.\par

  Notice that in the case when $ \psi(\cdot) \in \Psi(A,\infty)  $ and a function
 $ p \to p \cdot \log \psi(p) $ is convex,  then the space
$ G\psi $ coincides with some {\it exponential} Orlicz space. \par
 Conversely, if $ B < \infty, $ then the space $ G\psi(A,B) $ does  not coincides with
 the classical rearrangement invariant spaces: Orlicz, Lorentz, Marcinkiewicz  etc.\par

  The fundamental function of these spaces $ \phi(G(\psi), \delta) = ||I_A ||G(\psi), \mes(A) = \delta, \ \delta > 0, $
where  $ I_A  $ denotes as ordinary the indicator function of the measurable set $ A, $ by the formulae

$$
\phi(G(\psi), \delta) = \sup_{ p \in \supp (\psi)} \left[ \frac{\delta^{1/p}}{\psi(p)} \right]. \eqno(4.2)
$$

 The fundamental function of arbitrary rearrangement invariant spaces plays very important role in functional analysis,
theory of Fourier series and transform \cite{Bennet1} as well as in our further narration. \par

 Many examples of fundamental functions for some $ G\psi $ spaces are calculated in  \cite{Ostrovsky1}, \cite{Ostrovsky2}.\par

\vspace{3mm}

{\bf Remark 4.1} If we introduce the {\it discontinuous} function

$$
\psi_{(r)}(p) = 1, \ p = r; \psi_{(r)}(p) = \infty, \ p \ne r, \ p,r \in (A,B)
$$
and define formally  $ C/\infty = 0, \ C = \const \in R^1, $ then  the norm
in the space $ G(\psi_r) $ coincides with the $ L_r $ norm:

$$
||f||G(\psi_{(r)}) = |f|_r.
$$
Thus, the Grand Lebesgue Spaces are direct generalization of the
classical exponential Orlicz's spaces and Lebesgue spaces $ L_r. $ \par

\vspace{3mm}

{\bf Remark 4.2}  The function $ \psi(\cdot) $ may be generated as follows. Let $ \xi = \xi(x)$
be some measurable function: $ \xi: X \to R $ such that $ \exists  (A,B):
1 \le A < B \le \infty, \ \forall p \in (A,B) \ |\xi|_p < \infty. $ Then we can
choose

$$
\psi(p) = \psi_{\xi}(p) = |\xi|_p.
$$

 Analogously let $ \xi(t,\cdot) = \xi(t,x), t \in T, \ T $ is arbitrary set,
be some {\it family } $ F = \{ \xi(t, \cdot) \} $ of the measurable functions:
$ \forall t \in T  \ \xi(t,\cdot): X \to R $ such that
$$
 \exists  (A,B): 1 \le A < B \le \infty, \ \sup_{t \in T} \
|\xi(t, \cdot)|_p < \infty.
$$
 Then we can choose

$$
\psi(p) = \psi_{F}(p) = \sup_{t \in T}|\xi(t,\cdot)|_p. \eqno(4.3)
$$

 The function $ \psi_F(p) $ may be called as a {\it natural function} for the family $ F. $
This method was used in the probability theory, more exactly, in
the theory of random fields, see \cite{Kozachenko1},\cite{Ostrovsky1}, chapters 3,4. \par

 For instance, the function $ \Phi(\cdot) $  may be introduced by a natural way based on the
family

$$
F_{d,X} = \{ d(\xi(x_1),\xi(x_2)) \}, \ x_1,x_2 \in X.
$$

\vspace{3mm}

{\bf Remark 4.3} Note that the so-called {\it exponential} Orlicz spaces are particular cases of
Grand Lebesgue spaces  \cite{Kozachenko1}, \cite{Ostrovsky1}, p. 34-37.  In detail, let the $ N- $
Young-Orlicz function has a view

$$
 N(u) = e^{\mu(u)},
$$
where the function $ u \to \mu(u) $ is convex even twice differentiable function such that

$$
\lim_{u \to \infty} \mu'(u) = \infty.
$$
 Introduce a new function
$$
\psi_{\{N\}}(x) = \exp \left\{ \frac{\left[\log N(e^x) \right]^*}{x}   \right\},
$$
where $  g^*(\cdot) $ denotes the Young-Fenchel transform of the function $  g:  $

$$
g^*(x) = \sup_y (xy - g(y)).
$$
 Conversely,  the $  N  - $ function may be calculated up to equivalence
  through corresponding function $ \psi(\cdot) $  as follows:

 $$
 N(u) = e^{\tilde{\psi}^*(\log |u|) }, \ |u| > 3; \ N(u) = C u^2, |u| \le 3; \  \tilde{\psi}(p) = p \log \psi(p). \eqno(4.4)
 $$
 The Orlicz's space $ L(N) $ over our probabilistic space is equivalent up to sublinear norms equality with
Grand Lebesgue space $ G\psi_{\{N\}}. $ \par

\vspace{3mm}

{\bf Remark 4.4.} The theory of probabilistic {\it exponential} Grand Lebesgue spaces
or equally exponential Orlicz spaces gives a  very convenient apparatus for investigation of
the r.v. with exponential decreasing tails of distributions. Namely, the non-zero  r.v. $ \eta $ belongs to the
Orlicz space $ L(N), $  where $ N = N(u) $ is function described before, if and only if

$$
{\bf P} (\max(\eta, -\eta) > z) \le \exp(-\mu(C z)), \ z > 1,  \ C = C(N(\cdot), ||\eta||L(N)) \in (0,\infty).
$$
(Orlicz's version). \par

  Analogously may be written a Grand Lebesgue version of this inequality.
 In detail,  if $ 0 < ||\eta||G\psi< \infty,  $ then

 $$
{\bf P} (\max(\eta, -\eta) > z) \le  2\exp \left(- \tilde{\psi}(\log [z /||\eta||G\psi] ) \right), z \ge ||\eta||G\psi.
 $$
  Conversely, if

 $$
{\bf P} (\max(\eta, -\eta) > z) \le  2\exp \left(- \tilde{\psi}(\log [z /K] ) \right), z \ge K,
$$
then $ ||\eta||G\psi \le C(\psi)  \cdot K, \ C(\psi) \in (0,\infty). $  \par

\vspace{3mm}

 A very important subclass of the $ G\psi $ spaces form the so-called $ B(\phi) $ spaces. \par

Let $ \phi = \phi(\lambda), \lambda \in (-\lambda_0, \lambda_0), \ \lambda_0 = \const \in (0, \infty] $ be some even strong
convex which takes positive values for positive arguments twice continuous differentiable function, such that
$$
 \phi(0) = 0, \ \phi^{//}(0) \in(0,\infty), \ \lim_{\lambda \to \lambda_0} \phi(\lambda)/\lambda = \infty. \eqno(4.5)
$$
 We denote the set of all these function as $ \Phi; \ \Phi =\{ \phi(\cdot) \}. $ \par
 We say that the {\it centered} random variable (r.v) $ \xi  $
belongs to the space $ B(\phi), $ if there exists some non-negative constant
$ \tau \ge 0 $ such that

$$
\forall \lambda \in (-\lambda_0, \lambda_0) \ \Rightarrow
{\bf E} \exp(\lambda \xi) \le \exp[ \phi(\lambda \ \tau) ]. \eqno(4.6)
$$
 The minimal value $ \tau $ satisfying (4.6) is called a $ B(\phi) \ $ norm
of the variable $ \xi, $ write
 $$
 ||\xi||B(\phi) = \inf \{ \tau, \ \tau > 0: \ \forall \lambda \ \Rightarrow
 {\bf E}\exp(\lambda \xi) \le \exp(\phi(\lambda \ \tau)) \}.
 $$
 This spaces are very convenient for the investigation of the r.v. having a
exponential decreasing tail of distribution, for instance, for investigation of the limit theorem,
the exponential bounds of distribution for sums of random variables,
non-asymptotical properties, problem of continuous of random fields,
study of Central Limit Theorem in the Banach space etc.\par

  The space $ B(\phi) $ with respect to the norm $ || \cdot ||B(\phi) $ and
ordinary operations is a Banach space which is isomorphic to the subspace
consisted on all the centered variables of Orlicz's space $ (\Omega,F,{\bf P}), N(\cdot) $ with $ N \ - $ function

$$
N(u) = \exp(\phi^*(u)) - 1, \ \phi^*(u) = \sup_{\lambda} (\lambda u -\phi(\lambda)). \eqno(4.7)
$$
 The transform $ \phi \to \phi^* $ is called Young-Fenchel transform. The proof of considered assertion used the
properties of saddle-point method and theorem of Fenchel-Moraux:
$$
\phi^{**} = \phi.
$$

 The next facts about the $ B(\phi) $ spaces are proved in \cite{Kozachenko1}, \cite{Ostrovsky1}, p. 19-40:

$$
{\bf 1.} \ \xi \in B(\phi) \Leftrightarrow {\bf E } \xi = 0, \ {\bf and} \ \exists C = \const > 0,
$$

$$
U(\xi,x) \le \exp(-\phi^*(Cx)), x \ge 0,
$$
where $ U(\xi,x)$ denotes in this article the {\it tail} of distribution of the r.v. $ \xi: $

$$
U(\xi,x) = \max \left( {\bf P}(\xi > x), \ {\bf P}(\xi < - x) \right),
\ x \ge 0,
$$
and this estimation is  asymptotically exact. \par
 Here and further $ C, C_j, C(i) $ will denote the non-essentially positive
finite "constructive" constants.\par

 The function $ \phi(\cdot) $ may be "constructively" introduced by the formula
$$
\phi(\lambda) = \phi_0(\lambda) \stackrel{def}{=} \log \sup_{t \in T}
 {\bf E} \exp(\lambda \xi(t)), \eqno(4.8)
$$
 if obviously the family of the centered r.v. $ \{ \xi(t), \ t \in T \} $ satisfies the {\it uniform } Kramer's condition:
$$
\exists \mu \in (0, \infty), \ \sup_{t \in T} U(\xi(t), \ x) \le \exp(-\mu \ x), \ x \ge 0. \eqno(4.9)
$$
 In this case we will call the function $ \phi(\lambda) = \phi_0(\lambda) $ {\it natural } function. \par

 {\bf 2.} We define $ \psi(p) = \psi_{\phi}(p) := p/\phi^{-1}(p), \ p \ge 2. $
   It is proved that the spaces $ B(\phi) $ and $ G(\psi) $ coincides:$ B(\phi) = G(\psi) $ (set equality) and both
the norm $ ||\cdot||B(\phi) $ and $ ||\cdot|| $ are equivalent: $ \exists C_1 =
C_1(\phi), C_2 = C_2(\phi) = \const \in (0,\infty), \ \forall \xi \in B(\phi) $

$$
||\xi||G(\psi) \le C_1 \ ||\xi||B(\phi) \le C_2 \ ||\xi||G(\psi). \eqno(4.10)
$$

   The Gaussian (more precisely, subgaussian) case is considered in \cite{Garsia1}, \cite{Hu1}, \cite{Ral'chenko1}
 may be obtained by choosing $ \Phi(z) = \Phi_2(z) := \exp(z^2/2) - 1 $ or equally $ \psi(p) = \psi_2(p) = \sqrt{p}. $
 It may be considered easily more general example when
 $ \Phi(z) =  \Phi_Q(z) := \exp(|z|^Q/Q) - 1, \ Q = \const > 0; \ \Leftrightarrow \psi(p) = \psi_Q(p) := p^{1/Q}, \ p \ge 1. $  \par

 In the last case the following implication holds:

 $$
 \eta \in L(\Phi_Q), \ Q > 1 \ \Leftrightarrow  U(\eta,x) \le \exp \left( - C(\Phi, \eta) \ x^{Q'}  \right),
 $$
where as usually $ Q' = Q/(Q-1). $\par

\vspace{4mm}

 Assume that  the number $  \theta, $ measure $ m, $  distance $  d_{(\psi)},  $ and the function $ \psi = \psi(p) $  are such that
$  \theta > 0, \ C =  C(\theta) = \const \in (0,\infty); $

$$
(A,B):= \supp \psi(\cdot), \ A > \theta,  \ B > A;
$$

$$
 d_{(\psi)}(x_1,x_2) := ||\xi(x_1) - \xi(x_2))||G\psi;
$$

$$
 m^2(B(d_{(\psi)}, r,x)) \ge r^{\theta}/C(\theta), \ r \in [0,\diam(X, d_{(\psi)})], \  \ C(\theta) \in (0,\infty).
$$

 Define also a new   function:

 $$
\psi_{\theta}(p) \stackrel{def}{=}  \ \psi(p)/(1 - \theta/p), \ p \in (A, B).
 $$

 \vspace{4mm}

 {\bf Theorem 4.1.}   The sample paths of the r.f. $ \xi(x) $  belong a.e. for all the values $ p \in (A,B) $  to the H\"older space
 $ H(d_{(\psi)}^{1 - \theta/p} ). $ Moreover,

 $$
 \left| \sup_{d_{(\psi)}(x_1, x_2) > 0} \frac{|\xi(x_1) - \xi(x_2)|}{d_{(\psi)}^{1 - \theta/p}(x_1, x_2)} \right|_p \le
 C\cdot \psi_{\theta}(p), \ p \in (A,B). \eqno(4.11)
 $$

 As a consequence: let the semi - distance $ \rho = \rho(x_1,x_2) $ be in addition such that

 $$
 d_{(\psi)}^{ 1 - \theta/A }(\cdot, \cdot) << \rho(\cdot, \cdot), \eqno(4.12)
 $$
  then the r.f. $ \xi(\cdot) $ belongs  to the H\"older space  $  H^o(\rho) $  with probability one. \par

\vspace{4mm}

{\bf Proof.} We start from the relation (3.2):

$$
|\xi(x_1) - \xi(x_2)|  \le 12 \ Z^{1/p} \ 4^{1/p} \ C^{1/p}(\theta) \  \frac{d_p^{1-\theta/p}(x_1, x_2)}{1-\theta/p},
$$
where $ Z $ is a non - negative  r.v.  $ Z = Z_{(p)} $ is such that $ {\bf E} Z \le 1.  $

It follows from the direct definition of the norm in $ G\psi  $ spaces

$$
d_p(x_1, x_2) =   |\xi(x_1) - \xi(x_2)|_p \le \psi(p) \cdot d_{(\psi)}(x_1,x_2), \ p \in (A, B);
$$
and  we derive after substituting

$$
|\xi(x_1) - \xi(x_2)|  \le C_1(A,B) \cdot Z^{1/p} \cdot \psi(p) \cdot  \frac{d_{( \psi  )   }^{1-\theta/p}(x_1, x_2)}{1-\theta/p},
$$
or equally

$$
\left| \frac{|\xi(x_1) - \xi(x_2)|}{d_{(\psi)}^{1 - \theta/p}(x_1, x_2)} \right| \le  C_1(A,B) \cdot Z^{1/p} \cdot \frac{\psi(p)}{1 - \theta/p} =
C_1(A,B) \cdot Z^{1/p} \cdot \psi_{\theta}(p).
$$
 Since the right-hand of the last inequality does not dependent on the variables $  x_1, x_2, $

$$
\sup_{d_{(\psi)}(x_1, x_2) > 0}  \left| \frac{|\xi(x_1) - \xi(x_2)|}{d_{(\psi)}^{1 - \theta/p}(x_1, x_2)} \right| \le
C_1(A,B) \cdot Z^{1/p} \cdot \psi_{\theta}(p). \eqno(4.13)
$$

 It remains to calculate the $ L_p $ norm on both the sides of the last inequality. \par

\vspace{3mm}

{\bf Remark 4.5.}  The case when $ \psi(p) = \sqrt{p}  $ correspondents  to the Gaussian (more generally, subgaussian) random field $ \xi(x). $
The case $ \psi(p) = \exp(C p) $ appears in the articles \cite{Arnold1} and  \cite{Imkeller1}. However, in both these cases the
condition (2.5) is satisfied.\par

 In the case $ \psi(p) = \psi_{(r)}(p) $ we obtain the statement of theorem 3.1. as a  particular case. \par

\vspace{4mm}

 {\bf Theorem 4.2.}  Suppose that all the conditions of theorem 4.1 are satisfied. Then the r.f. $  \xi(x) $ satisfies
 the CLT in H\"older space $  H^o(\rho). $ \par

\vspace{4mm}

 {\bf Proof.} We apply the statement of theorem 4.1 to the r.f. $ S_n(\cdot). $\par

  Note first of all that

$$
d_{\psi_R}(x_1, x_2) :=   || S_n(x_1) - S_n(x_2)|| G\psi_R \le ||\xi(x_1) - \xi(x_2)|| = d_{\psi}(x_1, x_2). \eqno(4.14)
$$

 Therefore,

$$
m^2(B(d_{\psi_R}, r,x) \ge r^{\theta}/C(\theta), \ x \in X.
$$

Denote

$$
\psi_{\theta,R}(p) = \frac{\psi_R(p)}{1 - \theta/p}. \eqno(4.15)
$$
  It is known \cite{Liflyand1} that the function $ \psi_{\theta,R}(\cdot) $ belongs to the set $  \Psi = \{ \psi \} $
with  at the same support $ (A,B).  $ \par

We have by virtue of proposition of theorem 4.1

$$
\sup_n \left| \sup_{d_{(\psi_R)}(x_1, x_2) > 0} \frac{|S_n(x_1) - S_n(x_2)|}{d_{(\psi_R)}^{1 - \theta/p}(x_1, x_2)} \right|_p \le
 C\cdot \psi_{\theta,R}(p), \ p \in (A,B).
 $$

All the more so

$$
\sup_n \left| \sup_{d_{(\psi_R)}(x_1, x_2) > 0} \frac{|S_n(x_1) - S_n(x_2)|}{d_{(\psi_R)}^{1 - \theta/A}(x_1, x_2)} \right|_p \le
 C\cdot \psi_{\theta,R}(p), \ p \in (A,B), \eqno(4.16)
$$
and hence

$$
\sup_n \left| \left| \sup_{d_{(\psi_R)}(x_1, x_2) > 0} \frac{|S_n(x_1) - S_n(x_2)|}{d_{(\psi_R)}^{1 - \theta/A}(x_1, x_2)} \right| \right|G\psi_{\theta,R}
\le C < \infty. \eqno(4.17)
$$
 It remains to repeat the arguments using by the proof of theorem  3.1. \par

\vspace{4mm}

\section{ CLT in rectangle H\"older spaces via the fractional order Sobolev-Grand Lebesgue Spaces.}

\vspace{3mm}

 Let $ D $ be convex non-empty bounded closed domain with Lipschitz boundary in the whole space $ R^d, \ d = 1,2,\ldots, $
and let $ f: D \to R $ be measurable function. \par
 {\it We assume further for simplicity that} $ D = [0,1]^d. $ \par
 We denote and define $ |x| = (\sum_{i=1}^d x_i^2 )^{1/2},  \ \alpha  = \const \in (0,1], $

$$
|f|_p = |f|_{p,D} = \left[ \int_D |f(x)|^p \ dx  \right]^{1/p},  \ |u(\cdot, \cdot)|_p = |u(\cdot, \cdot)|_{p,D^2} =
$$

$$
 \left[ \int_D \int_D |u(x,y)|^p \ dx dy \right]^{1/p}, \ p = \const \ge 1,
$$

$$
  \omega(f,\delta) = \sup \{|f(x) - f(y)|: \ x,y \in D,  |x-y| \le \delta  \},  \ \delta  \in [0, \diam(D)], \eqno(5.0)
$$

$$
G_{\alpha}[f] (x,y)= \frac{f(x) - f(y) }{|x-y|^{\alpha}},  \hspace{5mm}  \nu(dx,dy) = \frac{dx dy}{|x-y|}, \eqno(5.1)
$$

$$
|u(\cdot, \cdot)|_{p,\nu} = |u(\cdot, \cdot)|_{p, \nu, D^2} = \left[ \int_D \int_D |u(x,y)|^p \ \nu(dx, dy) \right]^{1/p}, \eqno(5.2)
$$

$$
||f||W(\alpha,p) = | G_{\alpha}[f] (\cdot, \cdot) |_{p, \nu, D^2}. \eqno(5.3)
$$
 The norm  $ ||\cdot||W(\alpha,p), $ more precisely, semi-norm  is said to be {\it fractional } Sobolev's norm or similar
{\it Aronszajn, Gagliardo or Slobodeckij } norm; see, e.g. \cite{Nezzaa1}.\par

  If in the definition (5.3) instead the $ L_p(D^2)  $ stands another norm $ ||\cdot|| V(D^2), $ for instance,  Lorentz,
 Marcinkiewicz or Grand Lebesgue, (we recall its definition further), we obtain correspondingly the definition of the
 fractional  $ ||\cdot|| V(D^2) $ norm. \par

 The inequality

 $$
 |f(t) - f(s)| \le  8 \cdot 4^{1/p} \cdot \left[\frac{\alpha + 1/p}{\alpha - 1/p}\right]   \cdot |t-s|^{\alpha - 1/p} \cdot
 ||f||W(\alpha,p), \eqno(5.4)
 $$
or equally

$$
\omega(f,\delta) \le 8 \cdot 4^{1/p} \cdot \left[\frac{\alpha + 1/p}{\alpha - 1/p}\right]   \cdot \delta^{\alpha - 1/p} \cdot
\left[ \int_D \int_D \frac{|f(x) - f(y)|^p \ dx dy }{|x-y|^{\alpha p + 1}}   \right]^{1/p}, \eqno(5.5)
$$
which is true in the case $  d = 1 $ (the multidimensional case will be consider further), $ p > 1/\alpha, $ is called {\it fractional}
Sobolev, or Aronszajn, Gagliardo, Slobodeckij inequality.\par
 More precisely, the inequality (5.4) implies that the function $ f $ may be redefined  on the set of measure zero as a
continuous function  for which (5.4) there holds.\par
 Another look on the inequality (5.4): it may be construed as an imbedding theorem from the Sobolev fractional space into the
space of (uniform) continuous functions on the set $  D. $ \par
 The proof of the our version of inequality (5.4) may be obtained immediately from an article \cite{Hu1}, which based in turn on the
famous Garsia-Rodemich-Rumsey inequality, see \cite{Garsia1}.\par
 There are many generalizations of fractional Sobolev's imbedding theorem: on the Sobolev-Orlicz's spaces \cite{Adams1}, p. 253-364, on the
so-called {\it integer} Sobolev-Grand Lebesgue spaces  \cite{Ostrovsky100}, on the Lorentz and Marcinkiewicz spaces etc.\par
 The applications of these inequalities in the theory of random processes is investigated in the article \cite{Ostrovsky101}.\par

 The predicate that $ x \in D $ imply  $ x = \vec{x} = (x_1,x_2, \ldots, x_d), \ 0 \le x_i \le 1. $\par

 We define as in \cite{Ral'chenko1}, \cite{Hu1} the {\it rectangle difference}  operator $ \Box[f](\vec{x}, \vec{y}) =
 \Box[f](x,y),  \ x,y \in D, \ f:D \to R $ as follows.

$$
\Delta^{(i)}[f](x,y) := f(x_1,x_2, \ldots, x_{i-1}, y_i, x_{i+1}, \ldots, x_d) - f(x_1,x_2, \ldots, x_{i-1}, x_i, x_{i+1}, \ldots, x_d),
$$
with obvious modification when $ i=1 $ or $ i=d; $
$$
\Box[f](x,y) \stackrel{def}{=} \left\{ \otimes _{i=1}^d \Delta^{(i)} \right\} [f](x,y).\eqno(5.6)
$$
 For instance, if $ d=2, $ then

$$
\Box[f](x,y) = f(y_1,y_2) - f(x_1,y_2) - f(y_1,x_2) + f(x_1,x_2).
$$

 If the function $  f: [0,1]^d \to R $ is $ d $ times continuous differentiable, then

 $$
 \Box[f](\vec{x},\vec{y}) = \int_{x_1}^{y_1}
 \int_{x_2}^{y_2} \ldots   \int_{x_d}^{y_d} \frac{ \partial^d f}{\partial x_1  \partial x_2  \ldots   \partial x_d } \ dx_1 dx_2 \ldots   dx_d .
 $$

 The {\it rectangle module of continuity  } $  \Omega(f, \vec{\delta} ) =  \Omega(f, \delta )  $  for the
(continuous  a.e.) function $ f $ and vector $ \vec{\delta} = \delta = ( \delta_1, \delta_2, \ldots, \delta_d)
\in [0,1]^d  $  may be defined as well as ordinary module of continuity $ \omega(f,\delta) $ as follows:

$$
\Omega(f, \vec{\delta} ) \stackrel{def}{=} \sup \{ |\Box[f](x,y)|, \ (x,y): |x_i - y_i| \le \delta_i, \ i = 1,2,\ldots,d \}.
$$

   Let $ \vec{\alpha} = \{  \alpha_k \}, \ \alpha_k \in (0,1], \ k=1,2,\ldots,d; \ p > p_0 \stackrel{def}{=}\max_k (1/\alpha_k),
  \ M = \card  \{i, \alpha_i = \min_k \alpha_k\}, \delta_i = |x_i - y_i|, \ \vec{\delta} = \{ \delta_i \}, i = 1,2,\ldots,d; $

$$
\vec{x}^{\vec{\alpha}} := \prod_{i=1}^d x_i^{\alpha_i}, \  \vec{\delta}^{\pm 1/p} :=  \left[\prod_{i=1}^d  \delta_i  \right]^{\pm 1/p},
$$

$$
G_{\vec{\alpha}}[f] (x,y) =  \frac{\Box[f](x,y) }{  |(\vec{x}- \vec{y})^{\vec{\alpha}} | },  \hspace{5mm}  \nu(dx,dy) = \frac{\vec{dx} \vec{ dy}}{|x-y|},
$$

$$
||f||W( \vec{\alpha},p) = | G_{\vec{\alpha}}[f] (\cdot, \cdot) |_{p, \nu, D^2}.
$$

 The norm  $ ||\cdot||W( \vec{\alpha},p), $ more precisely, semi-norm  is said to be {\it multidimensional fractional } Sobolev's norm or similar
{\it Aronszajn, Gagliardo or Slobodeckij } norm.  \par
 Define also  the following function

$$
\zeta_{\vec{\alpha}}(p):=||f||W(\vec{\alpha},p), \  (A,B):= \supp \ [\zeta_{\vec{\alpha}}(\cdot)]
$$
 and suppose $  1 \le A < B \le \infty.  $ \par

  Denote $ A (\vec{\alpha}) = \max(A, p_0) $ and suppose also $ A(\vec{\alpha}) < B; $

 $$
 Q_{\alpha,d}(p) = 8^d \cdot 4^{d/p} \cdot \prod_{k=1}^d \left[ \frac{ \alpha_k + 1/p}{\alpha_k - 1/p} \right].
 $$

 We define a new  psi-function $ \psi_{\alpha}(p)  $ as follows.

$$
\psi_{\vec{\alpha}}(p):= \zeta_{\vec{\alpha}}(p) \cdot Q_{\alpha,d}(p).
$$

Let $ \xi = \xi(x) $ be again random  field.
 We introduce the following  {\it natural} $ \Psi $ function: $  \theta_{\vec{\alpha}}(p) = $

$$
\theta_{\alpha}(p) =  Q_{\alpha,d}(p) \cdot
\left[ \int_D \int_D {\bf E} |G_{\vec{\alpha}}[\xi](x,y)|^p \nu(dx,dy) \right]^{1/p}, \eqno(5.7)
$$

$$
\alpha = \vec{\alpha} = \{ \alpha_1, \alpha_2, \ldots,\alpha_d \}, \ \alpha_k = \const > 0;
$$
and suppose the function $ \theta_{\alpha}(p) $ has non-trivial support such that

$$
A = \inf \supp \theta_{\alpha}(\cdot) \ge 1/\min_k\alpha_k, \ B = \sup \supp \theta_{\alpha} \in (A, \infty].
$$

\vspace{3mm}

 {\bf Theorem 5.1.}

Let $ \nu(p) = \nu_{\alpha}(p) $ be some function from the set $ \Psi(A,B) $ such that the function
 $ \gamma(p) =\gamma_{\alpha}(p) = \nu(p)/\theta_{\alpha}(p) $ belongs also to the set $ \Psi(A,B). $  Then

$$
|| \Omega(\xi,\delta) ||G\nu \le \delta^{\alpha} \cdot \phi(G\gamma,1 /\delta). \eqno(5.8)
$$

\vspace{4mm}

{\bf Proof.} We will use the following result from \cite{Ostrovsky101}:

$$
| \Omega(\xi, \delta) |_p \le \delta^{\alpha - 1/p} \ \theta_{\alpha}(p), \ p \in (A,B),\eqno(5.9)
$$
from which follows

$$
\frac{| \Omega(\xi, \delta)|_p}{\nu(p) \cdot \delta^{\alpha}} \le  \frac{ (1/\delta)^{1/p}}{\gamma(p)}.  \eqno(5.10)
$$

 It remains to take supremum over $ p; \ p \in (A,B)  $ from both the sides of the last inequality (5.10). \par

\vspace{4mm}

{\bf Theorem 5.2.} Denote

$$
\theta_{\alpha, R}(p) = \frac{C_R \ p}{e \cdot \ln p} \cdot \theta_{\alpha}(p). \eqno(5.11)
$$
 Let $ \nu(p) = \nu_{\alpha}(p) $ be some function from the set $ \Psi(A,B) $ such that the function
 $ \gamma_R(p) =\gamma_{\alpha,R}(p) = \nu(p)/\theta_{\alpha,R}(p) $ belongs also to the set $ \Psi(A,B). $  Then

$$
|| \Omega(S_n,\delta) )||G\nu \le \delta^{\alpha} \cdot \phi(G\gamma_R,1 /\delta). \eqno(5.12)
$$

\vspace{3mm}

{\bf Proof} is alike to one in theorem 5.1., in which we substitute  the r.f. $ S_n(\cdot) $ instead  the r.f. $ \xi(\cdot) $
and apply the Rosenthal's inequality.\par

\vspace{4mm}

{\bf  Definition 5.1 of the rectangle H\"older space.  } \par

\vspace{3mm}

 Let $ f: \ D \to R $ be continuous function and let $ \omega = \omega(\delta) = \omega(\vec{\delta}), \ 0 \le \delta_i \le 1 $ be
some non - trivial rectangle module of continuity, i.e. non - negative continuous monotonically increasing over each variable   $  \delta_i $
function such that

$$
\omega(\delta) = 0 \ \Leftrightarrow \exists i = 1,2, \ldots, d:  \ \delta_i = 0. \eqno(5.13)
$$

 Define the following rectangle H\"older's norm

$$
||f||H_r(\omega) \stackrel{def}{=} \sup_{x \in D} |f(x)| + \sup_{\delta > 0} \left[ \frac{\Omega(f,\delta)}{\omega(\delta)} \right], \eqno(5.14)
$$
and correspondingly the rectangle H\"older's space $ H_r(\omega) $ which consists on all the (continuous) functions $ f: D \to R $
with finite norm $ ||f||H_r(\omega). $\par

 This space is not separable, therefore we define the (closed) its subspace (separable component)
 $ H_r^o(\omega)$ consisting on all the function from the space  $ H_r(\omega) $ satisfying the additional condition

$$
\lim_{|\vec{\delta}| \to 0+}  \left[ \frac{\Omega(f,\delta)}{\omega(\delta)} \right] = 0,
$$
under at the same norm. \par
  It follows immediately from theorem 5.2 the following assertion. \par

\vspace{3mm}

{\bf Theorem 5.3.}  We retain all the notations and conditions of theorem 5.2.  Suppose that the module of continuity
$ \omega_0 = \omega_0(\delta) $  be such that

$$
\lim_{|\vec{\delta}| \to 0} \left\{ \frac{ \omega_0(\delta) }{   \delta^{\alpha} \cdot \phi(G\gamma_R,1 /\delta)} \right\}  = \infty. \eqno(5.15)
$$
 Then the r.f. $ \xi(x) $ satisfies the CLT in the rectangle H\"older's space $ H_r^o(\omega_0). $ \par

\vspace{4mm}

  Authors does not know another versions of the CLT in the rectangle H\"older's spaces.\par

\vspace{4mm}

\section{Reducing of degree.}

\vspace{3mm}

 Let $ X = [0, 1]^m, \  m = 2, 3, . . . .$  In the articles  \cite{Kwapien1}, \cite{Ral'chenko1}, \cite{Hu1}
is obtained under some additional conditions (condition 2.5 etc.)
 a multivariate generalization of famous Garsia-Rodemich-Rumsey
inequality \cite{Garsia1}. Roughly speaking, instead degree "2" in our inequalities  stands degree 1 and
coefficients dependent on the distance $ d. $ \par

 The ultimate (sharp) value of this degree in general case of arbitrary metric space $ (X, d) $
is now unknown; see also \cite{Arnold1}, \cite{Imkeller1}. \par

\vspace{3mm}

{\it We intend to generalize the statement of theorem 3.1 on the case when the Young (Young-Orlicz) function $ \Phi $
satisfies in addition the condition (2.5) and is twice continuous differentiable.} \par

\vspace{3mm}

 Some new notations.  As in the second section

$$
d = d_{\Phi} = d_{\Phi}(x_1,x_2) := || \xi(x_1) - \xi(x_2)||L(\Phi).
$$
Further,
$$
K(\Phi) := \sup_{x,y > 0} \left[ \frac{\Phi^{-1}(xy)}{\Phi^{-1}(x) + \Phi^{-1}(y)}\right] < \infty;  \hspace{7mm}
C(\Phi) := \frac{\Phi^{-1}(1)}{54 K^2(\Phi)}; \eqno(6.1)
$$

 Let us define also the following important distance function: $ \tau(x_1, x_2) = $
 $$
  \tau(x_1, x_2; \Phi) = \tau(x_1, x_2; \Phi, m ) = \tau(x_1, x_2; \Phi, m,d) =
 $$

$$
  \max \left\{ \int_0^{d(x_1, x_2)} \Phi^{-1} \left[ \frac{1}{m(B(r,x_1)) } \right] \ dr, \
 \int_0^{d(x_1,x_2)} \Phi^{-1}  \left[ \frac{1}{m(B(r,x_2))} \right] \ dr \right\}. \eqno(6.2)
$$

  If $ \forall X_1, x_2 \in X \ \Rightarrow \tau(x_1, x_2) < \infty, $ then the measure $  m(\cdot) $ is called
{\it  weakly majorizing;} in this case the function $  \tau = \tau(x_1, x_2) $ satisfies the triangle inequality
and other properties of the distance function \cite{Kwapien1}. \par

 Let the function $ \phi(\cdot) $ be natural function for the r.f. $ \xi(x), \ x \in X: $
$$
\phi(\lambda) = \log \sup_{x \in X} {\bf E} \exp(\lambda \xi(x)),
$$
 if obviously the family of the centered r.v. $ \{ \xi(x), \ x \in X \} $ satisfies the {\it uniform } Kramer's condition.
 Denote

$$
\overline{\phi}(\lambda) = \sup_{n=1,2,\ldots} [n \phi(\lambda/\sqrt{n})],  \hspace{7mm} \overline{\Phi}(u) = \exp \left( \overline{\phi}^*(u) \right) - 1.
\eqno(6.3)
$$

\vspace{4mm}

{\bf Theorem 6.1.} Let $  m(\cdot) $ be some weakly majorizing measure relative the Young function  $ \overline{\Phi}(\cdot). $
 Suppose also $  K( \overline{\Phi}) < \infty. $ Let also $ \theta = \theta(x_1, x_2) $ be arbitrary distance function such that
$ \tau_{ \overline{\Phi}} << \theta. $ \par
  Then the r.f. $ \xi(x) $ satisfies CLT in H\"older space $ H^o(\theta). $ \par

 \vspace{3mm}

 {\bf Proof.} The inequality

$$
\left| \left| 0.5 \ C(\Phi) \ \sup_{\tau(x_1, x_2) > 0}  \frac{|\xi(x_1) - \xi(x_2)|}{\tau(x_1,x_2)} \right| \right| L(\Phi)\le 1 \eqno(6.4)
$$
for the r.f. $ \xi(x) $  is in fact proved in \cite{Kwapien1}; see also \cite{Ostrovsky107}. As we knew, $ \Phi(u) = \exp(\phi^*(u)) - 1.  $ \par
 We apply the last inequality (6.4) for the random field $ S_n(\cdot). $ Let us estimate first of all the moment generating function for the
 r.v. $ S_n(x_1) - S_n(x_2). $ \par

 Recall preliminarily that if $ \{ \eta_i \}, \ i = 1,2,\ldots $  be a sequence of i., i.d. centered r.v. satisfying the Kramer's condition:

 $$
  \exists \lambda_0 \in (0, \infty],  \ \forall \lambda: \ |\lambda| < \lambda_0 \ \Rightarrow
   \phi(\lambda) := \log {\bf E} \exp( \lambda \eta_1  ) < \infty,
 $$
then

$$
\sup_n \log {\bf E} \exp \left( n^{-1/2} \sum_{i=1}^n \eta_i  \right) \le \overline{\phi}(\lambda), \ |\lambda| < \lambda_0. \eqno(6.5)
$$

 We can rewrite the inequality (6.5) taking into account the relation between the functions $ \phi $ and $ \Phi $ as follows

$$
\sup_n || n^{-1/2} \sum_{i=1}^n \eta_i ||L(\overline{\Phi})  \le C_1 \ ||\eta_1||L(\Phi). \eqno(6.6)
$$

 Therefore,

$$
\sup_n || S_n(x_1) - S_n(x_2) ||L(\overline{\Phi})  \le C_1 \ d_{\Phi}(x_1, x_2), \eqno(6.7)
$$
and we conclude by means of the estimate  (6.4)

$$
\sup_n \left| \left| C(\Phi) \
\sup_{\tau(x_1, x_2) > 0}  \frac{|S_n(x_1) - S_n(x_2)|}{\tau(x_1,x_2)} \right| \right| L(\overline{\Phi})\le 1. \eqno(6.8)
$$

 It remains to repeat the arguments  using by the proof of theorem  3.1. \par


\vspace{4mm}

\begin{thebibliography}{99}

\vspace{3mm}

\bibitem{Adams1}
{\sc Adams R.A.} {\it Sobolev Spaces.} Academic Press, (1978), New York, San Francisco, London.

\bibitem{Arnold1}
{\sc Arnold L. and Imkeller P.}
 {\it On the spatial asymptotic Behavior of stochastic Flows in
Euclidean Space. }
Stoch. Processes Appl., {\bf 62(1),}  (1996), 19-54.

\bibitem{Barlow1}
{\sc Barlow M.T. and Yor M.} {\it Semimartingale inequalities via the Garsia-Rodemich-Rumsey lemma, and
applications to local times.} J. Funct. Anal.; {\bf 49(2),} (1982), 198-229.

\bibitem{Bednorz1}
 {\sc Bednorz W.} (2006). {\it A theorem on Majorizing Measures.}
 Ann. Probab., {\bf 34}, 1771-1781. MR1825156.

\bibitem{Bednorz2}
{\sc Bednorz W.} {\it The majorizing measure approach to the sample boundedness.}
arXiv:1211.3898v1 [math.PR] 16 Nov 2012

\bibitem{Bednorz3}
{\sc Bednorz W.}  (2010), {\it Majorizing measures on metric spaces. } C.R. math. Acad. Sci.
Paris, (2010), 348, no. 1-2, 75-78, MR2586748

\bibitem{Bednorz4}
{\sc Bednorz W.}  {\it H\"older continuity of random processes.}
arXiv:math/0703545v1 [math.PR] 19 Mar 2007.

\bibitem{Bennet1}
{\sc Bennett C. and Sharpley R.} {\it Interpolation of operators.}  Orlando, Academic Press Inc.,1988.

\bibitem{Dudley1}
{\sc Dudley R.M.} {\it Uniform Central Limit Theorem.}  Cambridge University Press, 1999.

\bibitem{Fernique1}
 {\sc Fernique X.} (1975). {\it Regularite des trajectoires des
    function aleatiores gaussiennes.}  Ecole de Probablite de
    Saint-Flour, IV - 1974, Lecture Notes in Mathematic. {\bf 480}, 1-96,
    Springer Verlag, Berlin.

\bibitem{Fernique2}
{\sc Fernique X,} {\it Caracterisation de processus de trajectoires majores ou
continues.} Seminaire de Probabilit´s XII. Lecture Notes in Math. 649,
(1978), 691–706, Springer, Berlin.

\bibitem{Fernique3}
{\sc Fernique X.} {\it Regularite de fonctions aleatoires non gaussiennes.} Ecolee
de Ete de Probabilit´s de Saint-Flour XI-1981. Lecture Notes in Math.
976, (1983), 1–74, Springer, Berlin.

\bibitem{Frolov1}
{\sc Frolov A.S., Tchentzov N.N.} {\it On the calculation by the Monte-Carlo
method definite integrals depending on the parameters.} Journal of Computetional
Mathematics and Mathematical Physics, (1962), V. 2, Issue 4, p. 714-718
(in Russian).

\bibitem{Fiorenza3}
 {\sc Fiorenza A., and Karadzhov G.E.} {\it Grand and small Lebesgue spaces and
       their analogs.} Consiglio Nationale Delle Ricerche, Instituto per le
      Applicazioni del Calcoto Mauro Picone, Sezione di Napoli, Rapporto tecnico n.
      272/03, (2005).

\bibitem{Garsia1}
{\sc Garsia, A. M.; Rodemich, E.; and Rumsey, H., Jr.} {\it A real variable lemma and
the continuity of paths of some Gaussian processes.} Indiana Univ. Math. J. 20
(1970/1971), 565-578.

\bibitem{Grigorjeva1}
{\sc Grigorjeva M.L., Ostrovsky E.I.}  {\it Calculation of Integrals on discontinuous
Functions by means of depending trials method.} Journal of Computetional
Mathematics and Mathematical Physics, (1996), V. 36, Issue 12, p. 28-39 (in
Russian).

\bibitem{Heinkel1}
{\sc Heinkel B.} {\it  Measures majorantes et le theoreme de la limite centrale dans $ C(S). $ }
 Z. Wahrscheinlichkeitstheory. verw. Geb., (1977). {\bf 38}, 339-351.

\bibitem{Hu1}
{\sc Yaozhong Hu and Khoa Le}
{\it A multiparameter Garsia-Rodemich-Rumsey inequality and some applications.}
arXiv:1211.6809v1 [math.PR] 29 Nov 2012

\bibitem{Imkeller1}
{\sc Imkeller P. and  Scheutzov M.}
 {\it Stratonovich calculus with spatial parameters and  anticipative problem in
multiplicative ergodic theory.}
Ann. Probab. {\bf 27(1),} (1999), 109-129.

\bibitem{Iwaniec2}
 {\sc Iwaniec T., P. Koskela P., and Onninen J.} {\it Mapping of finite distortion:
   Monotonicity and Continuity.}  Invent. Math. 144 (2001), 507-531.

\bibitem{Kassman1}
{\sc Kassman M.} {\it A Note on Integral Inequalities and Embedding of Besov Spaces. } Journal of Inequalities in
Pure and Applied Mathematics. V.4 Issue 4, article 107,   (2003), 47-57.

\bibitem{Kozachenko1}
 {\sc Kozachenko Yu. V., Ostrovsky E.I.} (1985). {\it The Banach Spaces of
      random Variables of subgaussian type.} Theory of Probab. and Math.
      Stat. (in Russian). Kiev, KSU, {\bf 32}, 43-57.

\bibitem{Krasnoselsky1}
{\sc Krasnoselsky M.A., Rutizky Ya.B.} {\it Convex function and Orlicz spaces.}
GIFML, Moskow, 1958 (in Russian).

\bibitem{Kwapien1}
{\sc Kwapien S. and Rosinsky J.} {\it Sample H\"older continuity of stochastic processes and majorizing measures.}
 (2004).  Seminar on Stochastic Analysis, Random Fields
and Applications IV, Progr. in Probab. 58, 155–163. Birkh\"ouser, Basel.

\bibitem{Ledoux1}
 {\sc Ledoux M., Talagrand M.} (1991) {\it Probability in Banach Spaces.}
      Springer, Berlin, MR 1102015.

\bibitem{Liflyand1}
{\sc Liflyand E., Ostrovsky E., Sirota L.} {\it Structural Properties of Bilateral Grand Lebesgue Spaces.}
Turk. J. Math.; {\bf 34} (2010), 207-219.

\bibitem{Nezzaa1}
{\sc Nezzaa E.D., Palatuccia G., Valdinocia E.}
{\it Hitchhiker’s guide to the fractional Sobolev spaces.}
arXiv:1104.4345v3 [math.FA] 19 Nov 2011

\bibitem{Ostrovsky1}
{\sc  Ostrovsky E.I.} (1999). {\it Exponential estimations for random Fields and its
applications (in Russian).}  Moscow-Obninsk, OINPE.

\bibitem{Ostrovsky2}
{\sc  Ostrovsky E. and Sirota L.}
{\it Moment Banach spaces: theory and applications.}
HIAT Journal of Science and Engineering, {\bf C}, Volume 4, Issues 1-2,
pp. 233-262, (2007).

\bibitem{Ostrovsky100}
{\sc  Ostrovsky E. and Sirota L.} {\it Module of continuity for the functions
 belonging to the Sobolev-Grand Lebesgue Spaces.}
arXiv:1006.4177v1 [math.FA] 21 Jun 2010

\bibitem{Ostrovsky101}
 {\sc Ostrovsky E., Sirota L.}
{\it Continuity of Functions belonging to the fractional Order Sobolev's-Grand Lebesgue Spaces.}
arXiv:1301.0132v1 [math.FA] 1 Jan 2013

\bibitem{Ostrovsky102}
{\sc Ostrovsky E., Rogover E.} {\it Exact exponential Bounds for the random Field Maximum Distribution
via the Majorizing Measures (Generic Chaining.)}
arXiv:0802.0349v1 [math.PR] 4 Feb 2008

\bibitem{Ostrovsky103}
{\sc Ostrovsky E.I.} (2002). {\it Exact exponential Estimations for Random Field
Maximum Distribution.} Theory Probab. Appl. 45 v.3, 281-286.

\bibitem{Ostrovsky104}
 {\sc Ostrovsky E., Sirota L.} {\it A counterexample to a hypothesis of light tail of maximum distribution
 for continuous random processes with light finite-dimensional tails.  }
arXiv:1208.6281v1 [math.PR] 30 Aug 2012

\bibitem{Ostrovsky105}
{\sc Ostrovsky E., Rogover E.} {\it Non-asymptotic exponential bounds for
MLE deviation under minimal conditions via classical and generic chaining methods.}
arXiv:0903.4062v1 [math.PR] 24 Mar 2009

\bibitem{Ostrovsky106}
 {\sc Ostrovsky E., Sirota L.}
{\it Monte-Carlo method for multiple parametric integrals
calculation and solving of linear integral Fredholm equations
of a second kind, with confidence regions in uniform norm.}
arXiv:1101.5381v1 [math.FA] 27 Jan 2011

\bibitem{Ostrovsky107}
 {\sc Ostrovsky E., Sirota L.} {\it Simplification of the majorizing method, with development.}
arXiv:1302.3202v1 [math.PR] 13 Feb 2013

\bibitem{Ostrovsky601}
{\sc Ostrovsky E.,  Sirota L.  } {\it Schl\"omilch and Bell series for Bessel's functions, with
probabilistic applications. }
arXiv:0804.0089v1 [math.CV] 1 Apr 2008


\bibitem{Ostrovsky602}
{\sc Ostrovsky E.,  Sirota L.  } {\it Uniform measures on the compact metric spaces, with applications.}
arXiv:1403.5725v1 [math.FA] 23 Mar 2014

\bibitem{Pisier1}
{\sc Pisier G.}  { Condition d'entropic assupant la continuite de certains processus et applications a l'analyse harmonique.  }
Seminaire d' analyse fonctionalle, (1980), Exp. 13, p. 22-33.

\bibitem{Prokhorov1}
{\sc Prokhorov Yu.V.}   {\it Convergence of Random Processes and Limit Theorems
of Probability Theory.} Probab. Theory Appl., (1956), V. 1, 177-238.

\bibitem{Ral'chenko1}
{\sc Ral'chenko, K. V. } {\it The two-parameter Garsia-Rodemich-Rumsey inequality and its
application to fractional Brownian fields.} Theory Probab. Math. Statist. No. 75
(2007), 167-178.

\bibitem{Runst1}
{\sc Runst T., Sickel W.}
{\it New Sobolev Spaces of Fractional Order, Nemytskij Operators, and Nonlinear Partial Didderential Equations.}
(1996), De Gruyter Incorporated, Walter;    Berlin, Heidelberg, London, New York, Hong Kong.

\bibitem{Talagrand1}
 {\sc Talagrand M.} (1996). {\it Majorizing measure: The generic chaining.}
 Ann. Probab., {\bf 24} 1049-1103. MR1825156

\bibitem{Talagrand2}
 {\sc Talagrand M.} (2005). {\it The Generic Chaining. Upper and
     Lower Bounds of Stochastic Processes.} Springer, Berlin. MR2133757.

\bibitem{Talagrand3}
{\sc Talagrand M.} (1987). Regularity of Gaussian processes. Acta Math. 159 no. 1-2, 99–
149, MR 0906527.

\bibitem{Talagrand4}
{\sc Talagrand M.} (1990), Sample boundedness of stochastic processes under increment
conditions. Annals of Probability 18, N. 1, 1-49, MR1043935.

\bibitem{Talagrand5}
{\sc Talagrand M.} (1992). A simple proof of the majorizing measure theorem. Geom.
Funct. Anal. 2, no. 1, 118-125. MR 1143666

\vspace{11mm}

\bibitem{Gusejnov1}
{\sc Gusejnov A.I., Muchtarov Ch.Sh.} {\it Introduction to the theory of non-linear singular integral equations.  }
 Moskow, Nauka, (1980), (in Russian).

\bibitem{Klicnarov'a1}
{\sc Klicnarov'a Jana.} {\it Central limit theorem for H\"older processes on $ R^m $ cube.    }
Comment.Math.Univ.Carolin. {\bf 48}, 1, (2007), 83-91.

\bibitem{Ratchkauskas1}
{\sc Ratchkauskas A, Suquet Ch.}
{\it Central limit theorems in H\"ölder topologies for Banach space valued random fields.}
Teor. Veroyatnost. i Primenen., 2004, Volume 49, Issue 1, Pages 109 \ - \ 125, (in Russian).

\bibitem{Ratchkauskas2}
{\sc A. Ratchkauskas, Ch. Suquet.} {\it Necessary and sufficient condition for the H\"olderian functional central limit theorem. }
J. Theoret. Probab. 17 (2004) 221–243.

\bibitem{Ratchkauskas3}
{\sc Ratchkauskas A, Suquet Ch.} {\it H\"older norm test statistics for epidemic change. } J. Statist. Plann. Inference, {\bf 126}, (2004),
495  \ - \ 520.

\bibitem{Ratchkauskas4}
{\sc Ratchkauskas A, Suquet Ch.} {\it Testing epidemic changes of infinite dimensional parameters.} Stat. Inference Stoch.
Process,  {\bf 9},  (2006), 111-134.

\bibitem{Ratchkauskas5}
{\sc Ratchkauskas A., V. Zemlys V.}
{\it Functional central limit theorem for a double-indexed summation process,} Liet. Mat.
Rink., {\bf 45},  (2005), 401-412.

\bibitem{Rosenthal1}
{\sc Rosenthal H.P.} {\it On the subspaces of $ L_p \ (p \ge 2) $ spanned by sequences of
independent Variables.} Israel J. Math., 1970, V.3 pp. 273-253.






\end{thebibliography}
\end{document}